\documentclass[11pt, a4paper]{article}
\usepackage{amsmath}
\usepackage{wasysym}
\usepackage{latexsym}
\usepackage{amsfonts}
\usepackage{mathrsfs}
\usepackage{amssymb}
\usepackage{ifsym}
\usepackage{dsfont}
\usepackage[all]{xy}
\usepackage{authblk}

\usepackage{siunitx} % Required for alignment

\newtheorem{Definitions1}{Definition}[section]

\newtheorem{Theorems1}{Theorem}[section]

\newtheorem{Coroll1}[Theorems1]{Corollary}
\newtheorem{Lemma1}[Theorems1]{Lemma}

\newtheorem{Quest1}{Question}[section]

\newenvironment{proof}[1][Proof]{\begin{trivlist}
\item[\hskip \labelsep {\bfseries #1}]}{\end{trivlist}}

\newcommand{\qed}{\nobreak \ifvmode \relax \else
      \ifdim\lastskip<1.5em \hskip-\lastskip
      \hskip1.5em plus0em minus0.5em \fi \nobreak
      \vrule height0.75em width0.5em depth0.25em\fi}

\begin{document}
\title{Partially-elementary end extensions of countable admissible sets}
\author{Zachiri McKenzie\\
Department of Philosophy, Zhejiang University\\
{\tt zach.mckenzie@gmail.com}}
\maketitle

\begin{abstract}
A result of Kaufmann \cite{kau81} shows that if $L_\alpha$ is countable, admissible and satisfies $\Pi_n\textsf{-Collection}$, then $\langle L_\alpha, \in \rangle$ has a proper $\Sigma_{n+1}$-elementary end extension. This paper investigates to what extent the theory that holds in $\langle L_\alpha, \in \rangle$ can be transferred to the partially-elementary end extensions guaranteed by Kaufmann's result. We show that there are $L_\alpha$ satisfying full separation, powerset and $\Pi_n\textsf{-Collection}$ that have no proper $\Sigma_{n+1}$-elementary end extension satisfying either $\Pi_{n}\textsf{-Collection}$ or $\Pi_{n+3}\textsf{-Foundation}$. In contrast, we show that if $A$ is a countable admissible set that satisfies $\Pi_n\textsf{-Collection}$ and $T$ is a recursively enumerable theory that holds in $\langle A, \in \rangle$, then $\langle A, \in \rangle$ has a proper $\Sigma_n$-elementary end extension that satisfies $T$. 
\end{abstract}

\section[Introduction]{Introduction}

In \cite{km68}, Keisler and Morley prove that every countable model of $\mathsf{ZF}$ has proper elementary end extension. Kaufmann \cite{kau81} refines this result showing that if $n \geq 1$ and $\mathcal{M}$ is a countable structure in the language of set theory that satisfies $\mathsf{KP}+\Pi_n\textsf{-Collection}$, then $\mathcal{M}$ has proper $\Sigma_{n+1}$-elementary end extension. And, conversely, if $n \geq 1$ and $\mathcal{M}$ is a structure in the language of set theory that satisfies $\mathsf{KP}+\mathsf{V=L}$ and has a proper $\Sigma_{n+1}$-elementary end extension, then $\mathcal{M}$ satisfies $\Pi_n\textsf{-Collection}$. In particular this shows that for a countable limit ordinal $\alpha$ and $n \geq 1$, $\langle L_\alpha, \in \rangle$ has a proper $\Sigma_{n+1}$-elementary end extension if and only if it satisfies $\Pi_n\textsf{-Collection}$. In the context of first-order arithmetic, the McDowell-Specker Theorem \cite{ms61} reveals that {\it every} model of $\mathsf{PA}$ has a proper elementary end extension. This is refined in Paris and Kirby \cite{pk78} where it is shown that if $n \geq 2$ and $\mathcal{M}$ is a countable structure in the language of arithmetic that satisfies $\mathsf{I}\Delta_0$, then $\mathcal{M}$ satisfies the arithmetic collection scheme for $\Sigma_n$-formulae if and only if $\mathcal{M}$ has a proper $\Sigma_n$-elementary end extension.   

A natural question to ask is how much of the theory of $\mathcal{M}$ satisfying $\mathsf{KP}+\Pi_{n}\textsf{-Collection}$ can be made to hold in a proper $\Sigma_{n+1}$-elementary end extension whose existence is guaranteed by Kaufmann's result? In particular, is there a proper $\Sigma_{n+1}$-elementary end extension of $\mathcal{M}$ that also satisfies $\mathsf{KP}+\Pi_{n}\textsf{-Collection}$? Or, if $\mathcal{M}$ is transitive, is there a proper $\Sigma_{n+1}$-elementary end extension of $\mathcal{M}$ that satisfies full induction for all set-theoretic formulae? In section 3 we show that the answers to the latter two of these questions is ``no". For $n \geq 1$, there is an $L_\alpha$ satisfying full separation, powerset and $\Pi_n\textsf{-Collection}$ that has no proper $\Sigma_{n+1}$-elementary end extension satisfying either $\Pi_n\textsf{-Collection}$ or $\Pi_{n+3}\textsf{-Foundation}$. A key ingredient is a generalisation of a result due to Simpson (see \cite[Remark 2]{kau81}) showing that if $n \geq 1$ and $\mathcal{M}$ is a structure in the language of set theory satisfying $\mathsf{KP}+\mathsf{V=L}$ that has $\Sigma_n$-elementary end extension satisfying enough set theory and with a new ordinal but no least new ordinal, then $\mathcal{M}$ satisfies $\Pi_n\textsf{-Collection}$. Here ``enough set theory" is either $\mathsf{KP}+\Pi_{n-1}\textsf{-Collection}$ or $\mathsf{KP}+\Pi_{n+2}\textsf{-Foundation}$. In section 4 we obtain a strong converse to this generalisation of Simpson's result for countable admissible sets using the Barwise Compactness Theorem. We show that if $A$ is a countable admissible set that satisfies $\Pi_n\textsf{-Collection}$ and $T$ is a recursively enumerable theory that holds in $\langle A, \in \rangle$, then $\langle A, \in \rangle$ has a $\Sigma_n$-elementary end extension that satisfies $T$ with a new ordinal but no least new ordinal. 

\section[Background]{Background}

Let $\mathcal{L}$ be the language of set theory-- the language whose only non-logical symbol is the binary relation $\in$. Let $\Gamma$ be a collection of $\mathcal{L}$-formulae.
\begin{itemize}
\item $\Gamma\textsf{-Separation}$ is the scheme that consists of the sentences 
$$\forall \vec{z} \forall w \exists y \forall x (x \in y \iff (x \in w \land \phi(x, \vec{z})),$$
for all formulae $\phi(x, \vec{z})$ in $\Gamma$. $\textsf{Separation}$ is the scheme that consists of these sentences for every formula $\phi(x, \vec{z})$ in $\mathcal{L}$.
\item $\Gamma\textsf{-Collection}$ is the scheme that consists of the sentences
$$\forall \vec{z} \forall w ((\forall x \in w) \exists y \phi(x, y, \vec{z}) \Rightarrow \exists c (\forall x \in w) (\exists y \in c)\phi(x, y, \vec{z})),$$
for all formulae $\phi(x, y, \vec{z})$ in $\Gamma$. $\textsf{Collection}$ is the scheme that consists of these sentences for every formula $\phi(x, y, \vec{z})$ in $\mathcal{L}$.
\item $\Gamma\textsf{-Foundation}$ is the scheme that consists of the sentences 
$$\forall \vec{z} (\exists x \phi(x, \vec{z}) \Rightarrow \exists y (\phi(y, \vec{z}) \land (\forall w \in y) \neg \phi(w, \vec{z}))),$$
for all formulae $\phi(x, \vec{z})$ in $\Gamma$. If $\Gamma= \{x \in z\}$, then the resulting axiom is referred to as $\textsf{Set-Foundation}$. $\textsf{Foundation}$ is the scheme that consists of these sentences for every formula $\phi(x, \vec{z})$ in $\mathcal{L}$.  
\end{itemize}
Let $T$ be a theory in a language that includes $\mathcal{L}$. Let $\Gamma$ be a class of $\mathcal{L}$-formulae. A formula is $\Gamma$ in $T$ or $\Gamma^T$ if it is provably equivalent in $T$ to a formula in $\Gamma$. A formula is $\Delta_n$ in $T$ or $\Delta_n^T$ if it is both $\Sigma_n^T$ and $\Pi_n^T$. 
\begin{itemize}
\item $\Delta_n\textsf{-Separation}$ is the scheme that consists of the sentences 
$$\forall \vec{z}(\forall v (\phi(v, \vec{z}) \iff \psi(v, \vec{z})) \Rightarrow \forall w \exists y \forall x(x \in y \iff (x \in w \land \phi(x, \vec{z}))))$$
for all $\Sigma_n$-formulae $\phi(x, \vec{z})$ and $\Pi_n$-formulae $\psi(x, \vec{z})$.
\item $\Delta_n\textsf{-Foundation}$ is the scheme that consists of the sentences 
$$\forall \vec{z}(\forall v(\phi(x, \vec{z}) \iff \psi(x, \vec{z})) \Rightarrow (\exists x \phi(x, \vec{z}) \Rightarrow \exists y (\phi(y, \vec{z}) \land (\forall w \in y) \neg \phi(w, \vec{z}))))$$
for all $\Sigma_n$-formulae $\phi(x, \vec{z})$ and $\Pi_n$-formulae $\psi(x, \vec{z})$. 
\end{itemize} 
Following \cite{mat01}, we take Kripke-Platek Set Theory ($\mathsf{KP}$) to be the $\mathcal{L}$-theory axiomatised by: \textsf{Extensionality}, \textsf{Emptyset}, \textsf{Pair}, \textsf{Union}, $\Delta_0\textsf{-Separation}$, $\Delta_0\textsf{-Collection}$ and $\Pi_1\textsf{-Foundation}$. Note that this differs from \cite{bar75, fri73}, which defines Kripke-Platek Set Theory to include \textsf{Foundation}. The theory $\mathsf{KPI}$ is obtained from $\mathsf{KP}$ by adding the axiom \textsf{Infinity}, which states that a superset of the von Neumann ordinal $\omega$ exists. We use $\mathsf{M}^-$ to denote the theory that is obtained from $\mathsf{KP}$ by replacing $\Pi_1\textsf{-Foundation}$ with $\textsf{Set-Foundation}$ and removing $\Delta_0\textsf{-Collection}$, and adding an axiom $\mathsf{TCo}$ asserting that every set is contained in a transitive set. The theory $\mathsf{M}$ is obtained from $\mathsf{M}^-$ by adding $\textsf{Powerset}$. Zermelo Set Theory ($\mathsf{Z}$) is obtained for $\mathsf{M}$ by removing $\mathsf{TCo}$ and adding $\textsf{Separation}$.

The theory $\mathsf{KP}$ proves $\mathsf{TCo}$ (see, for example, \cite[I.6.1]{bar75}). The following are some important consequences of fragments of the collection scheme over the theory $\mathsf{M}^-$:
\begin{itemize}
\item The proof of \cite[I.4.4]{bar75} generalises to show that, in the theory $\mathsf{M}^-$, $\Pi_n\textsf{-Collection}$ implies $\Sigma_{n+1}\textsf{-Collection}$.
\item \cite[Lemma 4.13]{flw16} shows that, over $\mathsf{M}^-$, $\Pi_n\textsf{-Collection}$ implies $\Delta_{n+1}\textsf{-Separation}$.
\item It is noted in \cite[Proposition 2.4]{flw16} that if $T$ is $\mathsf{M}^-+\Pi_n\textsf{-Collection}$, then the classes $\Sigma_{n+1}^T$ and $\Pi_{n+1}^T$ are closed under bounded quantification. 
\end{itemize} 

Let $\mathcal{M}= \langle M, \in^\mathcal{M} \rangle$ be an $\mathcal{L}$-structure. If $a \in M$, then we will use $a^*$ to denote the set $\{x \in M \mid \mathcal{M} \models (x \in a)\}$, as long as $\mathcal{M}$ is clear from the context. Let $\Gamma$ be a collection of $\mathcal{L}$-formulae. We say $X \subseteq M$ is $\Gamma$ over $\mathcal{M}$ if there is a formula $\phi(x, \vec{z})$ in $\Gamma$ and $\vec{a} \in M$ such that $\{x \in M \mid \mathcal{M} \models \phi(x, \vec{a})\}$. In the special can that $\Gamma$ is all $\mathcal{L}$-formulae, we say that $X$ is a {\bf definable subclass} of $\mathcal{M}$. A set $X \subseteq M$ is $\Delta_n$ over $\mathcal{M}$ if it is both $\Sigma_n$ over $\mathcal{M}$ and $\Pi_n$ over $\mathcal{M}$. 

A structure $\mathcal{N}= \langle N, \in^\mathcal{N} \rangle$ is an {\bf end extension} of $\mathcal{M}= \langle M, \in^\mathcal{M} \rangle$, written $\mathcal{M} \subseteq_e \mathcal{N}$, if $\mathcal{M}$ is a substructure of $\mathcal{N}$ and for all $x \in M$ and for all $y \in N$, if $\mathcal{N} \models (y \in x)$, then $y \in M$. An end extension $\mathcal{N}$ of $\mathcal{M}$ is {\bf proper} if $M \neq N$. We say that $\mathcal{N}$ is a {\bf $\Sigma_n$-elementary end extension} of $\mathcal{M}$, and write $\mathcal{M} \prec_{e, n} \mathcal{N}$, if $\mathcal{M} \subseteq_e \mathcal{N}$ and $\Sigma_n$ properties are preserved between $\mathcal{M}$ and $\mathcal{N}$.     

As shown in \cite[Chapter II]{bar75}, the theory $\mathsf{KP}$ is capable of defining G\"{o}del's constructible universe ($L$). For all sets $X$,
$$\mathrm{Def}(X)= \{Y \subseteq X \mid Y \textrm{ is a definable subclass of } \langle X, \in \rangle \},$$
which can be seen to be a set in the theory $\mathsf{KP}$ using a formula for satisfaction in set structures such as the one described in \cite[Section III.1]{bar75}. The levels of $L$ are then defined by the recursion:
$$L_0= \emptyset \textrm{ and } L_\alpha= \bigcup_{\beta < \alpha} L_\beta \textrm{ if } \alpha \textrm{ is a limit ordinal,}$$
$$L_{\alpha+1}= L_\alpha \cup \mathrm{Def}(L_\alpha), \textrm{ and}$$
$$L= \bigcup_{\alpha \in \mathrm{Ord}} L_\alpha.$$
The function $\alpha \mapsto L_\alpha$ is total and $\Delta_1^{\mathsf{KP}}$. The axiom $\mathsf{V=L}$ asserts that every set is the member of some $L_\alpha$. A transitive set $M$ such that $\langle M, \in \rangle$ satisfies $\mathsf{KP}$ is said to be an {\bf admissible set}. An ordinal $\alpha$ is said to be an {\bf admissible ordinal} if $L_\alpha$ is an admissible set. 

Let $T$ be an $\mathcal{L}$-theory. A transitive set $M$ is said to be a {\bf minimum model} of $T$ if $\langle M, \in \rangle \models T$ and for all transitive sets $N$ with $\langle N, \in \rangle \models T$, $M \subseteq N$. For example, $L_{\omega_1^\mathrm{ck}}$ is the minimum model of $\mathsf{KPI}$. Gostanian \cite[\S 1]{gos80} shows that all sufficiently strong subsystems of $\mathsf{ZF}$ and $\mathsf{ZF}^-$ obtained by restricting the separation and collection schemes to formulae in the L\'{e}vy classes have minimum models. In particular:

\begin{Theorems1}
(Gostanian \cite{gos80}) Let $n \in \omega$. The theory $\mathsf{Z}+\Pi_n\textsf{-Collection}$ has a minimum model. Moreover, the minimum model of this theory satisfies $\mathsf{V=L}$. 
\end{Theorems1}   

The fact that $\mathsf{KP}$ is able to define satisfaction in set structures also facilitates the definition of formulae expressing satisfaction, in the universe, for formulae in any given level of the L\'{e}vy hierarchy.

\begin{Definitions1} \label{Df:Delta0Satisfaction}
The formula $\mathrm{Sat}_{\Delta_0}(q, x)$ is defined as
$$\begin{array}{c}
(q \in \omega) \land (q= \ulcorner \phi(v_1, \ldots, v_m) \urcorner \textrm{ where } \phi \textrm{ is } \Delta_0) \land (x= \langle x_1, \ldots, x_m \rangle) \land\\
\exists N \left( \bigcup N \subseteq N \land (x_1, \ldots, x_m \in N) \land (\langle N, \in \rangle \models \phi[x_1, \ldots, x_m]) \right)
\end{array}.$$
\end{Definitions1}

We can now inductively define formulae $\mathrm{Sat}_{\Sigma_n}(q, x)$ and $\mathrm{Sat}_{\Pi_n}(q, x)$ that express satisfaction for formulae in the classes $\Sigma_n$ and $\Pi_n$.

\begin{Definitions1}
The formulae $\mathrm{Sat}_{\Sigma_n}(q, x)$ and $\mathrm{Sat}_{\Pi_n}(q, x)$ are defined recursively for $n>0$. $\mathrm{Sat}_{\Sigma_{n+1}}(q, x)$ is defined as the  formula
$$\exists \vec{y} \exists k \exists b \left( \begin{array}{c}
(q= \ulcorner\exists \vec{u} \phi(\vec{u}, v_1, \ldots, v_l)\urcorner \textrm{ where } \phi \textrm{ is } \Pi_n)\land (x= \langle x_1, \ldots, x_l \rangle)\\
\land (b= \langle \vec{y}, x_1, \ldots, x_l \rangle) \land (k= \ulcorner \phi(\vec{u}, v_1, \ldots, v_l) \urcorner) \land \mathrm{Sat}_{\Pi_n}(k, b)
\end{array}\right);$$
and  $\mathrm{Sat}_{\Pi_{n+1}}(q, x)$ is defined as the formula
$$\forall \vec{y} \forall k \forall b \left( \begin{array}{c}
(q= \ulcorner\forall \vec{u} \phi(\vec{u}, v_1, \ldots, v_l) \urcorner \textrm{ where } \phi \textrm{ is } \Sigma_n)\land (x= \langle x_1, \ldots, x_l \rangle)\\
\land ((b= \langle \vec{y}, x_1, \ldots, x_l \rangle) \land (k= \ulcorner\phi(\vec{u}, v_1, \ldots, v_l)\urcorner) \Rightarrow \mathrm{Sat}_{\Sigma_n}(k, b))
\end{array}\right).$$
\end{Definitions1}

\begin{Theorems1} \label{Complexityofpartialsat} 
Suppose $n \in \omega$ and $m=\max \{ 1, n \}$. The formula $\mathrm{Sat}_{\Sigma_n}(q, x)$ (respectively $\mathrm{Sat}_{\Pi_n}(q, x)$) is $\Sigma_m^{\mathsf{KP}}$ ($\Pi_m^{\mathsf{KP}}$, respectively). Moreover, $\mathrm{Sat}_{\Sigma_n}(q, x)$ (respectively $\mathrm{Sat}_{\Pi_n}(q, x)$) expresses satisfaction for $\Sigma_n$-formulae ($\Pi_n$-formulae, respectively) in the theory $\mathsf{KP}$, i.e., if $\mathcal{M} \models\mathsf{KP}$, $\phi(v_1,\ldots,v_k)$ is a $\Sigma_n$-formula, and $x_1,\ldots,x_k$ are in $M$, then for $q = \ulcorner   \phi( v_1, \ldots, v_k) \urcorner$, $\mathcal{M}$ satisfies the universal generalisation of the following formula:
$$x= \langle x_1, \ldots,x_k \rangle \Rightarrow \left( \phi(x_1,\ldots,x_k) \iff \mathrm{Sat}_{\Sigma_{n}}(q, x) \right).$$
\Square
\end{Theorems1}

Friedman \cite[Section 2]{fri73} classifies the countable ordinals that can appear as the order type of the ordinals of a standard part of a nonstandard model of $\mathsf{KP}$. The key ingredient in Friedman's classification is the fact that every countable admissible set has an end extension with no least new ordinal that satisfies $\mathsf{KP}$.

\begin{Theorems1} \label{Th:FriedmanEndExtensionResult}
(Friedman \cite[Theorem 2.2]{fri73}) Let $M$ be a countable admissible set. Let $T$ be a recursively enumerable $\mathcal{L}$-theory such that $\langle M, \in \rangle \models T$. Then there exists $\mathcal{N}= \langle N, \in^\mathcal{N} \rangle$ such that $\langle M, \in \rangle \subseteq_e \mathcal{N}$, $\mathcal{N} \models T$ and $\mathrm{Ord}^{\mathcal{N}} \backslash \mathrm{Ord}^{\langle M, \in \rangle}$ is nonempty and has no least element. \Square
\end{Theorems1}

Barwise \cite[Appendix]{bar75} introduces the machinery of admissible covers to apply infinitary compactness arguments, such as the one used in the proof of Theorem \ref{Th:FriedmanEndExtensionResult}, to nonstandard countable models. The proof of \cite[Theorem A.4.1]{bar75} shows that for any countable model $\mathcal{M}$ of $\mathsf{KP}+\textsf{Foundation}$ and for any recursively enumerable $\mathcal{L}$-theory $T$ that holds in $\mathcal{M}$, $\mathcal{M}$ has proper end extension that satisfies $T$. By calibrating \cite[Appendix]{bar75}, Ressayre \cite[Theorem 2.15]{res87} shows that this result also holds for countable models of $\mathsf{KP}+\Sigma_1\textsf{-Foundation}$. 

\begin{Theorems1}
Let $\mathcal{M}= \langle M, \in^\mathcal{M} \rangle$ be a countable model of $\mathsf{KP}+\Sigma_1\textsf{-Foundation}$. Let $T$ be a recursively enumerable theory such that $\mathcal{M} \models T$. Then there exists $\mathcal{N} \models T$ such that $\mathcal{M} \subseteq_e \mathcal{N}$ and $M \neq N$. \Square 
\end{Theorems1}

Kaufmann \cite{kau81} identifies necessary and sufficient conditions for models of $\mathsf{M}^-$ to have proper $\Sigma_n$-elementary end extensions.

\begin{Theorems1} \label{Th:KaufmannEEResult} 
(Kaufmann \cite[Theorem 1]{kau81}) Let $n \geq 1$. Let $\mathcal{M}= \langle M, \in^\mathcal{M} \rangle$ be a model of $\mathsf{KP}$. Consider
\begin{itemize}
\item[(I)] there exists $\mathcal{N}= \langle N, \in^\mathcal{N} \rangle$ such that $\mathcal{M} \prec_{e, n+1} \mathcal{N}$ and $M \neq N$;
\item[(II)] $\mathcal{M} \models \Pi_{n}\textsf{-Collection}$.  
\end{itemize}
If $\mathcal{M} \models \mathsf{V=L}$, then $(I) \Rightarrow (II)$. If $M$ is countable, then $(II) \Rightarrow (I)$. \Square   
\end{Theorems1}

It should be noted that Kaufmann proves that (II) implies (I) in the above under the weaker assumption that $\mathcal{M}$ is a resolvable model of $\mathsf{M}^-$. A model $\mathcal{M}=\langle M, \in^\mathcal{M} \rangle$ of $\mathsf{M}^-$ is {\bf resolvable} if there is a function $F$ that is $\Delta_1$ over $\mathcal{M}$ such that for all $x \in M$, there exists $\alpha \in \mathrm{Ord}^\mathcal{M}$ such that $x \in F(\alpha)$. The function $\alpha \mapsto L_\alpha$ witnesses the fact that any model of $\mathsf{KP}+\mathsf{V=L}$ is resolvable.

\section[Admissible sets admitting topless partially elementary end extensions]{Admissible sets admitting topless partially elementary end extensions}

The next result is a generalisation of the result, due to Simpson, that that is mentioned in \cite[Remark 2]{kau81}. The proof of this generalisation is based on Enayat's proof of a refinement of Simpson's result (personal communication) that corresponds to the specific case of the following theorem when $n=1$ and $\mathcal{M}$ is transitive.

\begin{Theorems1} \label{Th:CollectionFromPartiallyElementaryEE1}
Let $n \geq 1$. Let $\mathcal{M}= \langle M, \in^\mathcal{M} \rangle$ be a model of $\mathsf{KP}+\mathsf{V=L}$. Suppose $\mathcal{N}= \langle N, \in^\mathcal{N} \rangle$ is such that $\mathcal{M} \prec_{e, n} \mathcal{N}$, $\mathcal{N} \models \mathsf{KP}$ and $\mathrm{Ord}^\mathcal{N} \backslash \mathrm{Ord}^\mathcal{M}$ is nonempty and has no least element. If $\mathcal{N}\models \Pi_{n-1}\textsf{-Collection}$ or $\mathcal{N} \models \Pi_{n+2}\textsf{-Foundation}$, then $\mathcal{M} \models \Pi_n\textsf{-Collection}$.
\end{Theorems1}

\begin{proof}
Assume that $\mathcal{N}= \langle N, \in^\mathcal{N} \rangle$ is such that 
\begin{itemize}
\item[(I)] $\mathcal{M} \prec_{e, n} \mathcal{N}$;
\item[(II)] $\mathcal{N} \models \mathsf{KP}$;
\item[(III)] $\mathrm{Ord}^\mathcal{N} \backslash \mathrm{Ord}^\mathcal{M}$ is nonempty and has no least element. 
\end{itemize}
Note that, since $\mathcal{M} \prec_{e, 1} \mathcal{N}$ and $\mathcal{M} \models \mathsf{V=L}$, for all $\beta \in \mathrm{Ord}^\mathcal{N} \backslash \mathrm{Ord}^\mathcal{M}$, $M \subseteq (L_\beta^{\mathcal{N}})^*$. We need to show that if either $\Pi_{n-1}\textsf{-Collection}$ or $\Pi_{n+2}\textsf{-Foundation}$ hold in $\mathcal{N}$, then $\mathcal{M} \models \Pi_n\textsf{-Collection}$. Let $\phi(x, y, \vec{z})$ be a $\Pi_n$-formula. Let $\vec{a}, b \in M$ be such that 
$$\mathcal{M} \models (\forall x \in b) \exists y \phi(x, y, \vec{a}).$$
So, for all $x \in b^*$, there exists $y \in M$ such that
$$\mathcal{M} \models \phi(x, y, \vec{a}).$$
Therefore, since $\mathcal{M} \prec_{e, n} \mathcal{N}$, for all $x \in b^*$, there exists $y \in M$ such that
$$\mathcal{N} \models \phi(x, y, \vec{a}).$$
Now, $\phi(x, y, \vec{z})$ can be written as $\forall w \psi(w, x, y, \vec{z})$ where $\psi(w, x, y, \vec{z})$ is $\Sigma_{n-1}$. Let $\xi \in \mathrm{Ord}^\mathcal{N} \backslash \mathrm{Ord}^\mathcal{M}$. So, for all $\beta \in \mathrm{Ord}^\mathcal{N} \backslash \mathrm{Ord}^\mathcal{M}$ and for all $x \in b^*$, there exists $y \in (L_\beta^\mathcal{N})^*$ such that
$$\mathcal{N} \models (\forall w \in L_\xi)\psi(w, x, y, \vec{a}).$$
Therefore, for all $\beta \in \mathrm{Ord}^\mathcal{N} \backslash \mathrm{Ord}^\mathcal{M}$,
\begin{equation}\label{eq:Equation1}
\mathcal{N} \models (\forall x \in b)(\exists y \in L_\beta)(\forall w \in L_\xi) \psi(w, x, y, \vec{a})
\end{equation}
Now, define $\theta(\beta, \xi, b, \vec{a})$ to be the formula
$$(\forall x \in b)(\exists y \in L_\beta)(\forall w \in L_\xi) \psi(w, x, y, \vec{a}).$$
If $\Pi_{n-1}\textsf{-Collection}$ holds in $\mathcal{N}$, then $\theta(\beta, \xi, b, \vec{a})$ is equivalent to a $\Sigma_{n-1}$-formula. Without $\Pi_{n-1}\textsf{-Collection}$, $\theta(\beta, \xi, b, \vec{a})$ can be written as a $\Pi_{n+2}$-formula. Therefore, $\Pi_{n-1}\textsf{-Collection}$ or $\Pi_{n+2}\textsf{-Foundation}$ in $\mathcal{N}$ will ensure that there is a least $\beta_0 \in \mathrm{Ord}^\mathcal{N}$ such that $\mathcal{N} \models \theta(\beta_0, \xi, b, \vec{a})$. Moreover, by (\ref{eq:Equation1}), $\beta_0 \in M$. Therefore,
$$\mathcal{N} \models (\forall x \in b) (\exists y \in L_{\beta_0})(\forall w \in L_\xi) \psi(w, x, y, \vec{a}).$$
So, for all $x \in b^*$, there exists $y \in (L_{\beta_0}^\mathcal{M})^*$, for all $w \in (L_\xi^\mathcal{N})^*$, 
$$\mathcal{N} \models \psi(w, x, y, \vec{a}).$$
Which, since $\mathcal{M} \prec_{e, n} \mathcal{N}$, implies that for all $x \in b^*$, there exists $y \in (L_{\beta_0}^\mathcal{M})^*$, for all $w \in M$, 
$$\mathcal{M} \models \psi(w, x, y, \vec{a}).$$
Therefore, $\mathcal{M} \models (\forall x \in b)(\exists y \in L_{\beta_0}) \phi(x, y, \vec{a})$. This shows that $\Pi_n\textsf{-Collection}$ holds in $\mathcal{M}$.    
\Square
\end{proof}  

Enayat uses a specific case of Theorem \ref{Th:CollectionFromPartiallyElementaryEE1} to show that the \mbox{$\langle L_{\omega_1^{\mathrm{ck}}}, \in \rangle$} has no proper $\Sigma_1$-elementary end extension that satisfies $\mathsf{KP}$ (personal communication). We now turn generalising this result to show that for all $n \geq 1$, the minimum model of $\mathsf{Z}+\Pi_n\textsf{-Collection}$ has no proper $\Sigma_{n+1}$-elementary end extension that satisfies either $\mathsf{KP}+\Pi_{n+3}\textsf{-Foundation}$ or $\mathsf{KP}+\Pi_n\textsf{-Collection}$. However, by Theorem \ref{Th:KaufmannEEResult}, for all $n \geq 1$, the minimum model of $\mathsf{Z}+\Pi_n\textsf{-Collection}$ does have a proper $\Sigma_{n+1}$-elementary end extension.

The following result follows from \cite[Theorem 4.4]{mck19}:

\begin{Theorems1}
Let $n \geq 1$. The theory $\mathsf{M}+\Pi_{n+1}\textsf{-Collection}+\Pi_{n+2}\textsf{-Foundation}$ proves that there exists a transitive model of $\mathsf{Z}+\Pi_n\textsf{-Collection}$. \Square
\end{Theorems1}

\begin{Coroll1} \label{Th:LimitsOfCollectionInMinimalModels}
Let $n \geq 1$. Let $M$ be the minimal model of $\mathsf{Z}+\Pi_{n}\textsf{-Collection}$. Then there is an instance of $\Pi_{n+1}\textsf{-Collection}$ that fails in $\langle M, \in \rangle$.\Square
\end{Coroll1}

\begin{Theorems1}
Let $n \geq 1$. Let $M$ be the minimal model of $\mathsf{Z}+\Pi_{n}\textsf{-Collection}$. Then $\langle M, \in \rangle$ has a proper $\Sigma_{n+1}$-elementary end extension, but neither
\begin{itemize}
\item[(I)] a proper $\Sigma_{n+1}$-elementary end extension satisfying $\mathsf{KP}+\Pi_{n+3}\textsf{-Foundation}$, nor
\item[(II)] a proper $\Sigma_{n+1}$-elementary end extension satisfying $\mathsf{KP}+\Pi_n\textsf{-Collection}$.
\end{itemize}
\end{Theorems1}

\begin{proof}
The fact that $\langle M, \in \rangle$ has a proper $\Sigma_{n+1}$-elementary end extension follows from Theorem \ref{Th:KaufmannEEResult}. Let $\mathcal{N}= \langle N, \in^{\mathcal{N}} \rangle$ be such that $\mathcal{N} \models \mathsf{KP}$, $N \neq M$ and $\langle M, \in \rangle \prec_{e, n+1} \mathcal{N}$. Since $M$ is the minimal model of $\mathsf{Z}+\Pi_{n}\textsf{-Collection}$, $\langle M, \in \rangle \models \neg \sigma$ where $\sigma$ is the sentence
$$\exists x (x \textrm{ is transitive} \land \langle x, \in \rangle \models \mathsf{M}+\Sigma_{n+1}\textsf{-Separation}+\Pi_{n}\textsf{-Collection}).$$
Since $\sigma$ is $\Sigma_1^{\mathsf{KP}}$ and $\langle M, \in \rangle \prec_{e, 1} \mathcal{N}$, $\mathcal{N} \models \neg \sigma$. Since $\mathcal{N} \models \mathsf{KP}$ and $M \neq N$, $\mathrm{Ord}^\mathcal{N} \backslash \mathrm{Ord}^{\langle M, \in\rangle}$ is nonempty. If $\gamma$ is the least element of $\mathrm{Ord}^\mathcal{N} \backslash \mathrm{Ord}^{\langle M, \in\rangle}$, then 
$$\mathcal{N} \models (\langle L_{\gamma}, \in \rangle \models \mathsf{Z}+\Pi_{n}\textsf{-Collection}),$$
which contradicts the fact that $\mathcal{N} \models \neg \sigma$. Therefore, $\mathrm{Ord}^\mathcal{N} \backslash \mathrm{Ord}^{\langle M, \in\rangle}$ is nonempty and contains no least element. Therefore, by Theorem \ref{Th:CollectionFromPartiallyElementaryEE1} and Corollary \ref{Th:LimitsOfCollectionInMinimalModels}, there must be both an instance of $\Pi_n\textsf{-Collection}$ and an instance of $\Pi_{n+3}\textsf{-Foundation}$ that fails in $\mathcal{N}$.     
\Square
\end{proof}

\section[Building partially elementary end extensions]{Building partially elementary end extensions}

In this section we show that if $M$ is a countable admissible set that satisfies $\Pi_n\textsf{-Collection}$, then we can build proper $\Sigma_n$-elementary end extensions of $\langle M, \in \rangle$ that satisfy as much of the first-order theory $\langle M, \in \rangle$ as we want them to. More specifically, we show that if $M$ is a countable admissible set with $\langle M, \in \rangle \models \Pi_n\textsf{-Collection}$ and $T$ is recursively enumerable with $\langle M, \in \rangle \models T$, then there exists a proper $\Sigma_n$-elementary end extension of $\langle M, \in \rangle$ that satisfies $T$.

We construct proper $\Sigma_n$-elementary end extensions of admissible sets using an appropriate version of the Barwise Compactness Theorem. In order to present this construction, it is convenient to introduce a family of class theory extensions of $\mathsf{KP}$. Closely related class theory extensions of $\mathsf{KP}$ have been used in \cite{fri73} to present Barwise Compactness Arguments and \cite{js18} to study extensions of $\mathsf{KP}$ obtained by adding fixed point axioms. Let $\mathcal{L}^c$ be the first-order language of class theory-- the language obtained from $\mathcal{L}$ adding a unary relation $\mathcal{S}$ that distinguishes sets from classes. To simplify the presentation of $\mathcal{L}^c$-formulae, we will treat $\mathcal{L}^c$ as a two-sorted language with sorts {\it sets} (referred to using lower case Roman letters $w, x, y, z, \ldots$) that are the elements of the domain that satisfy $\mathcal{S}$ and {\it classes} (referred to using upper case Roman letters $W, X, Y, Z \ldots$) that are any element of the domain. Therefore, $\exists x(\cdots )$ is an abbreviation for $\exists x(\mathcal{S}(x) \land \cdots)$, $\forall x(\cdots)$ is an abbreviation for $\forall x(\mathcal{S}(x) \Rightarrow \cdots)$, $\exists x(x=X)$ is an abbreviation for $\mathcal{S}(X)$, etc. We say that an $\mathcal{L}^c$-formula, $\phi$, is {\it elementary} if $\phi$ contains only atomic formulae in the form $x \in Y$, $x \in y$ and $\mathcal{S}(Y)$ and all of the quantifiers in $\phi$ are restricted to sets. In other words, an elementary $\mathcal{L}^c$-formula is a formula that does not contain the symbol $=$ and only contains set variables with the possible exception of subformulae in the form $x \in Y$ where $Y$ is a free (class) variable. The collection $\Delta_0^c$ is the smallest class of elementary $\mathcal{L}^c$-formulae that contains all atomic formulae, is closed under the connectives of propositional logic, and quantification in the form $\forall x \in y$ and $\exists x \in y$ where $x$ and $y$ are distinct variables. The classes $\Sigma_n^c$ and $\Pi_n^c$ are the classes of elementary formulae defined inductively from the class $\Delta_0^c$ in the usual way.
\begin{itemize}
\item $\mathsf{KP}^c$ is the $\mathcal{L}^c$-theory with axioms:
\begin{itemize}
\item[] $\forall X \forall Y (X \in Y \Rightarrow \exists x(x=X))$;
\item[]($\textsf{Extensionality}^c$) $\forall X \forall Y (X=Y \iff \forall x (x \in X \iff x \in Y))$;
\item[]($\textsf{Pairing}^c$) $\forall x \forall y \exists z \forall w(w \in z \iff w= x \lor w= y)$;
\item[]($\textsf{Union}^c$) $\forall x \exists y \forall z(z \in y \iff (\exists w \in x)(z \in w))$;
\item[]($\Delta_0^c\textsf{-Separation}$) for all $\Delta_0^c$-formulae, $\phi(x, \vec{Z})$,
$$\forall \vec{Z} \forall w \exists y \forall x(x \in y \iff (x \in w) \land \phi(x, \vec{Z}));$$
\item[]($\Delta_0^c\textsf{-Collection}$) for all $\Delta_0^c$-formulae, $\phi(x, y, \vec{Z})$,
$$\forall \vec{Z} \forall w((\forall x \in w) \exists y \phi(x, y, \vec{Z}) \Rightarrow \exists c (\forall x \in w) (\exists y \in c) \phi(x, y, \vec{Z}));$$
\item[]($\Pi_1^c\textsf{-Foundation}$) for all $\Pi_1^c$-formulae, $\phi(x, \vec{Z})$,
$$\forall \vec{Z}(\exists x \phi(x, \vec{Z}) \Rightarrow \exists y(\phi(y, \vec{Z}) \land (\forall w \in y) \neg \phi(w, \vec{Z})));$$
\item[]($\Delta_1^c\textsf{-}\mathsf{CA}$) for all $\Sigma_1^c$-formulae, $\phi(x, \vec{Z})$, and for all $\Pi_1^c$-formulae, $\psi(x, \vec{W})$,
$$\forall \vec{Z} \forall \vec{W}(\forall x(\phi(x, \vec{Z}) \iff \psi(x, \vec{W})) \Rightarrow \exists X \forall y (y \in X \iff \phi(x, \vec{Z}))).$$ 
\end{itemize}
\end{itemize}    
Friedman's class theory $\mathrm{Adm}^c$ \cite[Definition 1.14]{fri73} differs from $\mathsf{KP}^c$ by including the full scheme of foundation for $\mathcal{L}$-formulae instead of $\Pi_1^c\textsf{-Foundation}$. It should also be noted that the theory $\mathsf{KP}^c$ utilised in \cite{js18} includes full $\in$-induction for all $\mathcal{L}^c$-formulae. We have chosen to include only $\Pi_1^c\textsf{-Foundation}$ in $\mathsf{KP}^c$ in order to ensure that $\mathsf{KP}^c$ is a conservative extension of $\mathsf{KP}$, which here only includes $\Pi_1\textsf{-Foundation}$.

The theory $\mathsf{KP}^c$ proves that the class of elementary $\mathcal{L}^c$-formulae that are equivalent to a $\Sigma_1^c$-formula is closed under quantification that is bounded by a set variable and, similarly, the class of elementary $\mathcal{L}^c$-formulae that are equivalent to a $\Pi_1^c$-formula is also closed under quantification that is bounded by a set variable. The usual argument showing that $\mathsf{KP}$ proves $\Sigma_1\textsf{-Collection}$ adapts to show that $\mathsf{KP}^c$ proves $\Sigma_1^c\textsf{-Collection}$.

We will also be interested in extensions of $\mathsf{KP}^c$ that are obtained by strengthening the class comprehension scheme. For $n > 1$, define:
\begin{itemize}
\item[]($\Delta_n^c\textsf{-}\mathsf{CA}$) for all $\Sigma_n^c$-formulae, $\phi(x, \vec{Z})$, and for all $\Pi_n^c$-formulae, $\psi(x, \vec{W})$, 
$$\forall \vec{Z} \forall \vec{W}(\forall x(\phi(x, \vec{Z}) \iff \psi(x, \vec{W})) \Rightarrow \exists X \forall y (y \in X \iff \phi(x, \vec{Z}))).$$
\end{itemize}

We now turn to showing that the sets of any model of $\mathsf{KP}^c+\Delta_{n+1}\textsf{-}\mathsf{CA}$ satisfy $\mathsf{KP}+\Pi_n\textsf{-Collection}+\Pi_{n+1}\textsf{-Foundation}$, and conversely every model of $\mathsf{KP}+\Pi_n\textsf{-Collection}+\Pi_{n+1}\textsf{-Foundation}$ can be expanded to a model of $\mathsf{KP}^c+\Delta_{n+1}\textsf{-}\mathsf{CA}$. In particular, $\mathsf{KP}^c+\Delta_{n+1}\textsf{-}\mathsf{CA}$ is a conservative extension of $\mathsf{KP}+\Pi_n\textsf{-Collection}+\Pi_{n+1}\textsf{-Foundation}$ for sentences in the language of set theory.

\begin{Definitions1}
We use $\Delta_0^*$ to denote the smallest class of $\mathcal{L}$-formulae that contains the atomic formulae in the form $x \in y$, is closed under the connectives of propositional logic, and quantification in the form $\forall x \in y$ and $\exists x \in y$ where $x$ and $y$ are distinct variables. 
\end{Definitions1}

Note that the class $\Delta_0^*$, when viewed as a class of $\mathcal{L}^c$-formulae, is just the class of $\Delta_0^c$-formulae in which all variables are restricted to sets.

\begin{Lemma1} \label{Th:EqualityEliminationForDelta0}
Let $\phi(\vec{z})$ be a $\Delta_0$-formula. There is a $\Delta_0^*$-formula $\phi^\prime(\vec{z})$ such that 
$$\mathsf{Extensionality} \vdash \forall \vec{z} (\phi(\vec{z}) \iff \phi^\prime(\vec{z})).$$
\end{Lemma1}

\begin{proof}
Replace any subformula in the form $x=y$ with $(\forall w \in x)(w \in y)\land (\forall w \in y)(w \in x)$. \Square
\end{proof}

\begin{Theorems1}
Let $n \in \omega$. Let $\mathcal{M}= \langle M, \in^\mathcal{M}, \mathcal{S}^\mathcal{M} \rangle$ be a model of $\mathsf{KP}^c+\Delta_{n+1}\textsf{-}\mathsf{CA}$. Then $\mathcal{M}_{\mathrm{Set}}= \langle \mathcal{S}^\mathcal{M}, \in^\mathcal{M} \rangle$ satisfies $\mathsf{KP}+\Pi_n\textsf{-Collection}+\Pi_{n+1}\textsf{-Foundation}$.  
\end{Theorems1}

\begin{proof}
Let $\mathcal{M}_{\mathrm{Set}}= \langle \mathcal{S}^\mathcal{M}, \in^\mathcal{M} \rangle$. It is immediate that $\mathcal{M}_{\mathrm{Set}}$ satisfies $\textsf{Extensionality}$, $\textsf{Emptyset}$, $\textsf{Pair}$ and $\textsf{Union}$. Lemma \ref{Th:EqualityEliminationForDelta0} and the scheme of $\Delta_0^c\textsf{-Separation}$ in $\mathcal{M}$ imply that $\mathcal{M}_{\mathrm{Set}}$ satisfies $\Delta_0\textsf{-Separation}$. Similarly, employing Lemma \ref{Th:EqualityEliminationForDelta0} shows that $\Delta_0^c\textsf{-Collection}$ in $\mathcal{M}$ implies $\Delta_0\textsf{-Collection}$ in $\mathcal{M}_{\mathrm{Set}}$, and $\Pi_1^c\textsf{-Foundation}$ in $\mathcal{M}$ implies $\Pi_1\textsf{-Foundation}$ in $\mathcal{M}_{\mathrm{Set}}$. This shows that the theorem holds when $n=0$. Therefore, assume that $n > 0$. We need to verify that $\Pi_n\textsf{-Collection}$ and $\Pi_{n+1}\textsf{-Foundation}$ hold in $\mathcal{M}_{\mathrm{Set}}$. Let $V_{\Pi_n} \in M$ be such that
$$\langle \ulcorner \phi(x) \urcorner, a \rangle \in V_{\Pi_n} \textrm{ if and only if } \mathcal{M}_{\mathrm{Set}} \models \mathrm{Sat}_{\Pi_n}(\ulcorner \phi(x) \urcorner, a).$$
Let $V_{\Sigma_n} \in M$ be such that
$$\langle \ulcorner \phi(x) \urcorner, a \rangle \in V_{\Sigma_n} \textrm{ if and only if } \mathcal{M}_{\mathrm{Set}} \models \mathrm{Sat}_{\Sigma_n}(\ulcorner \phi(x) \urcorner, a).$$
Note that $\Delta_{n+1}\textsf{-}\mathsf{CA}$ in $\mathcal{M}$ ensures that the classes $V_{\Pi_n}$ and $V_{\Sigma_n}$ exist. To see that $\mathcal{M}_{\mathrm{Set}}$ satisfies $\Pi_n\textsf{-Collection}$, let $\phi(x, y, \vec{z})$ be a $\Pi_n$-formula. Let $b, \vec{a} \in \mathcal{S}^\mathcal{M}$ be such that 
$$\mathcal{M}_{\mathrm{Set}} \models (\forall x \in b)\exists y \phi(x, y, \vec{a}).$$
Consider $\theta(x, y, \vec{a}, v, V_{\Pi_n})$ defined by:
$$\exists u(u= \langle x, y, \vec{a} \rangle \land \langle v, u \rangle \in V_{\Pi_n}).$$
Note that, by Lemma \ref{Th:EqualityEliminationForDelta0}, $\theta(x, y, \vec{a}, v, V_{\Pi_n})$ is equivalent to a $\Sigma_1^c$-formula. Moreover, 
$$\mathcal{M} \models (\forall x \in b) \exists y \theta(x, y, \vec{a}, \ulcorner \phi(x, y, \vec{z}) \urcorner, V_{\Pi_n}).$$
Therefore, by $\Sigma_1^c\textsf{-Collection}$ in $\mathcal{M}$, 
$$\mathcal{M} \models \exists c(\forall x \in b) (\exists y \in c) \theta(x, y, \vec{a}, \ulcorner \phi(x, y, \vec{z}) \urcorner, V_{\Pi_n}).$$
Therefore,
$$\mathcal{M}_{\mathrm{Set}} \models \exists c(\forall x \in b)(\exists y \in c) \phi(x, y, \vec{a}).$$
This shows that $\mathcal{M}_{\mathrm{Set}}$ satisfies $\Pi_n\textsf{-Collection}$. 

Finally, we need to verify that $\mathcal{M}_{\mathrm{Set}}$ satisfies $\Pi_{n+1}\textsf{-Foundation}$. To this end, let $\phi(x, \vec{z})$ be a $\Pi_{n+1}$-formula. Therefore $\phi(x, \vec{z})$ can be written as $\forall \vec{y} \psi(\vec{y}, x, \vec{z})$ where $\psi(\vec{y}, x, \vec{z})$ is $\Sigma_n$. Let $\vec{a} \in \mathcal{S}^\mathcal{M}$.
Consider $\theta(\vec{y}, x, \vec{a}, v, V_{\Sigma_n})$ defined by:
$$\forall u (u= \langle \vec{y}, x, \vec{a} \rangle \Rightarrow \langle v, u \rangle \in V_{\Sigma_n}).$$
Note that $\theta(\vec{y}, x, \vec{a}, v, V_{\Sigma_n})$ is equivalent to a $\Pi_1^c$-formula. Therefore, $\Pi_1^c\textsf{-Foundation}$, the class
$$\{ x \mid \forall \vec{y} \theta(\vec{y}, x, \vec{a}, \ulcorner \psi(\vec{y}, x, \vec{z}) \urcorner, V_{\Sigma_n})\}$$
is either empty or has an $\in^\mathcal{M}$-least element in $\mathcal{M}$. Therefore, the class $\{x \mid \phi(x, \vec{a}) \}$ is either empty or has an $\in^\mathcal{M}$-least element in $\mathcal{M}_{\mathrm{Set}}$. This shows that $\mathcal{M}_{\mathrm{Set}}$ satisfies $\Pi_{n+1}\textsf{-Foundation}$.  
\Square
\end{proof}

Conversely, if $\mathcal{M}$ is a model of $\mathsf{KP}+\Pi_n\textsf{-Collection}+\Pi_{n+1}\textsf{-Foundation}$, then one can adjoin the classes that are $\Delta_{n+1}$ over $\mathcal{M}$ to $\mathcal{M}$ to obtain a model of $\mathsf{KP}^c+\Delta_{n+1}\textsf{-}\mathsf{CA}$.

\begin{Theorems1} \label{Th:ExpansionSatisfyingClassTheory}
Let $n \in \omega$. Let $\mathcal{M}= \langle M, \in^\mathcal{M} \rangle$ be a model of $\mathrm{KP}+\Pi_n\textsf{-Collection}+\Pi_{n+1}\textsf{-Foundation}$. Let 
$$\mathcal{X}= \{ X \subseteq M \mid X \textrm{ is } \Delta_{n+1} \textrm{ over } \mathcal{M} \} \backslash \{a^* \mid a \in M\}$$ $$\textrm{and } \in^\prime= \in^\mathcal{M} \cup \in \upharpoonright (M \times \mathcal{X}).$$
Then $\langle M \cup \mathcal{X}, \in^\prime, M \rangle \models \mathrm{KP}^c+\Delta_{n+1}\textrm{-}\mathrm{CA}$. 
\end{Theorems1}

\begin{proof}
Note that $\langle M \cup \mathcal{X}, \in^\prime, M \rangle$ clearly satisfies $\forall X \forall Y (X \in Y \Rightarrow \exists x(x=X))$, $\textsf{Extensionality}^c$, $\textsf{Pairing}^c$ and $\textsf{Union}^c$. Let $\phi(\vec{x}, Z_0, \ldots, Z_{m-1})$ be a $\Delta_0^c$-formula and let \mbox{$A_0, \ldots, A_{m-1} \in M \cup \mathcal{X}$}. Since for all $i \in m$, the formula $y \in A_i$ can be expressed as a $\Delta_{n+1}$-formula with parameters from $M$, there exists a $\Sigma_{n+1}$-formula $\psi(\vec{x}, \vec{z})$ and a $\Pi_{n+1}$-formula $\theta(\vec{x}, \vec{z})$, and $\vec{a} \in M$ such that for all $\vec{x} \in M$,
$$\langle M \cup \mathcal{X}, \in^\prime, M \rangle \models \phi(\vec{x}, A_0, \ldots, A_{m-1})$$ 
$$\textrm{if and only if } \mathcal{M} \models \psi(\vec{x}, \vec{a}) \textrm{ if and only if } \mathcal{M} \models \theta(\vec{x}, \vec{a}).$$ 
Therefore, $\Delta_0^c\textsf{-Separation}$ in $\langle M \cup \mathcal{X}, \in^\prime, M \rangle$ follows from $\Delta_{n+1}\textsf{-Separation}$ in $\mathcal{M}$, $\Delta_0^c\textsf{-Collection}$ in $\langle M \cup \mathcal{X}, \in^\prime, M \rangle$ follows from $\Sigma_{n+1}\textsf{-Collection}$ in $\mathcal{M}$, and $\Pi_1^c\textsf{-Foundation}$ in $\langle M \cup \mathcal{X}, \in^\prime, M \rangle$ follows from $\Pi_{n+1}\textsf{-Foundation}$ in $\mathcal{M}$. Similarly, if $\phi(x, \vec{Z})$ is a $\Sigma_{n+1}^c$-formula ($\Pi_{n+1}^c$-formula) and $\vec{A} \in M \cup \mathcal{X}$, then there exists a $\Sigma_{n+1}$-formula ($\Pi_{n+1}$-formula, respectively), $\psi(x, \vec{z})$ and $\vec{a} \in M$ such that for all $x \in M$,
$$\langle M \cup \mathcal{X}, \in^\prime, M \rangle \models \phi(x, \vec{A}) \textrm{ if and only if } \mathcal{M} \models \psi(x, \vec{a}).$$
Therefore, $\langle M \cup \mathcal{X}, \in^\prime, M \rangle$ satisfies $\Delta_{n+1}\textsf{-}\mathsf{CA}$.
\Square
\end{proof} 

We will build $\Sigma_n$-elementary end extensions of countable admissible sets satisfying $\Pi_n\textsf{-Collection}$ using an appropriate version of the Barwise Compactness Theorem. 

\begin{Definitions1}
Let $A$ be a countable admissible set. Let $\mathcal{L}^\prime$ be obtained from $\mathcal{L}$ by adding constant symbols $\bar{a}$ for each $a \in A$. Let $\mathcal{L}^{\prime\prime}$ be obtained from $\mathcal{L}^\prime$ by adding constant symbols $\mathbf{c}_n$ for each $n \in \omega$. Write $\mathcal{L}^{\prime\prime}_A$ for the fragment of $\mathcal{L}^{\prime\prime}_{\omega_1 \omega}$ that is coded in $A$. We will identified an $\mathcal{L}^{\prime\prime}_A$-theory, $T$, with the set of codes of sentences in $T$. 
\end{Definitions1}

The following version of the Barwise Compactness Theorem appears as \cite[Theorem 1.12]{fri73}:

\begin{Theorems1} \label{Th:BarwiseCompactness}
(Barwise Compactness) Let $A$ be a countable admissible set and let $\mathcal{X} \subseteq \mathcal{P}(A)$ such that $\langle A \cup \mathcal{X}, \in, A \rangle \models \mathsf{KP}^c$. If $T$ is an $\mathcal{L}^{\prime\prime}_A$-theory with $T \in A \cup \mathcal{X}$ and for all $T_0 \subseteq T$ with $T_0 \in A$, $T_0$ has a model, then $T$ has a model. \Square 
\end{Theorems1}

The proof of the next result is based on the proof of \cite[Theorem 2.2]{fri73}.

\begin{Theorems1} \label{Th:MainEndExtensionResult}
Let $n \in \omega$. Let $A$ be a countable admissible set such that $\langle A, \in \rangle \models \Pi_n\textsf{-Collection}$. Let $S$ be a recursively enumerable $\mathcal{L}$-theory such that $\langle A, \in \rangle \models S$. Then there exists $\mathcal{M}= \langle M, \in^\mathcal{M} \rangle$ such that
\begin{itemize}
\item[(i)] $\mathcal{M} \models S$;
\item[(ii)] $\langle A, \in \rangle \prec_n \mathcal{M}$;
\item[(iii)] $\mathrm{Ord}^\mathcal{M} \backslash A$ has no least element. 
\end{itemize} 
\end{Theorems1}

\begin{proof}
Let
$$\mathcal{X}= \{ X \subseteq A \mid X \textrm{ is a } \Delta_{n+1} \textrm{ subset of } A \} \backslash A.$$
By Theorem \ref{Th:ExpansionSatisfyingClassTheory}, $\langle A \cup \mathcal{X}, \in, A \rangle \models \mathsf{KP}^c+\Delta_{n+1}\textsf{-}\mathsf{CA}$. Let $T_0$ be $\mathcal{L}_A^{\prime\prime}$-theory with axioms:
\begin{itemize}
\item[(I)] $S$;
\item[(II)] for all $a \in A$,
$$\forall x\left(x \in \bar{a} \iff \bigvee_{b \in a} x= \bar{b}\right);$$
\item[(III)] for all $\alpha \in \mathrm{Ord}^{\langle A, \in \rangle}$, $\bar{\alpha} \in \mathbf{c}_0$;
\item[(IV)] for all $\Pi_n$-formulae, $\phi(x_1, \ldots, x_m)$ and for all $a_1, \ldots, a_m \in A$ such that $\langle A, \in \rangle \models \phi(a_1, \ldots, a_m)$,
$$\phi(\bar{a_1}, \ldots \bar{a_m}).$$ 
\end{itemize}
Note that $T_0 \in A \cup \mathcal{X}$ and for all $T^\prime \subseteq T_0$ with $T^\prime \in A$, $T^\prime$ has a model. Therefore, by Theorem \ref{Th:BarwiseCompactness}, $T_0$ is consistent. Let $\langle \phi_n \mid k \in \omega \rangle$ be an enumeration of the $\mathcal{L}_A^{\prime\prime}$-sentences. Let $m \in \omega$ and suppose that $T_{3m} \supseteq T_0$ has been defined, contains only finitely many of the constant symbols $\mathbf{c}_k$ and is consistent. Define
$$T_{3m+1}= \left\{ \begin{array}{ll}
T_{3m} \cup \{\phi_m\} & \textrm{ if } T_{3m} \cup \{\phi_m\} \textrm{ is consistent}\\
T_{3m} \cup \{\neg \phi_m\} & \textrm{ otherwise}  
\end{array}\right.$$
Define:
$$\textrm{if } \neg \phi_m \in T_{3m+1} \textrm{ and } \phi_m= \bigwedge \Gamma, \textrm{ then let}$$
$$T_{3m+2}= T_{3m+1} \cup \{\neg \psi\} \textrm{ for some } \psi \in \Gamma \textrm{ with } T_{3m+1} \cup \{\neg \psi\} \textrm{ consistent;}$$
$$\textrm{if } \neg \phi_m \in T_{3m+1} \textrm{ and } \phi_m= \forall x \psi(x), \textrm{ then let}$$   
$$T_{3m+2}= T_{3m+1} \cup \{\neg \psi(\mathbf{c}_k)\} \textrm{ for some } \mathbf{c}_k \textrm{ not appearing in } T_{3m+1};$$
$$T_{3m+2}= T_{3m+1} \textrm{ otherwise.}$$
Define:
$$\textrm{if for some } \alpha \in \mathrm{Ord}^{\langle A, \in \rangle}, T_{3m+2}\cup \{\neg(\bar{\alpha} \in c_m)\} \textrm{ is consistent, then let}$$ 
$$T_{3m+3}= T_{3m+2} \cup \{\neg(\bar{\alpha} \in c_m)\};$$
$$\textrm{otherwise let } T_{3m+3}= T_{3m+2} \cup \{(\mathbf{c}_k \in \mathbf{c}_m)\}\cup \{ \bar{\alpha} \in \mathbf{c}_k \mid \alpha \in \mathrm{Ord}^{\langle A, \in \rangle}\}$$ 
$$\textrm{where } \mathbf{c}_k \textrm{ does not appear in } T_{3m+2}.$$
Now, $T_{3m+3} \supseteq T_{3m} \supseteq T_0$ and $T_{3m+3}$ only contains finitely many of the constant symbols $\mathbf{c}_k$. Moreover, by Theorem \ref{Th:BarwiseCompactness}, $T_{3m+3}$ is consistent. Now, let
$$T= \bigcup_{m \in \omega} T_m.$$
Therefore, $T$ is consistent. Now, the terms of $T$ form a Henkin model $\mathcal{M}=\langle M, \in^\mathcal{M} \rangle$ satisfying (i)-(iii).   
\Square
\end{proof}

Combined with Theorem \ref{Th:CollectionFromPartiallyElementaryEE1}, we obtain the following characterisation of admissible $L_\alpha$ that have proper $\Sigma_n$-elementary end extensions with no least new ordinal satisfying any recursively enumerable fragment of the theory of $\langle L_\alpha, \in \rangle$.

\begin{Coroll1}
Let $L_\alpha$ be countable and admissible. Then the following are equivalent:
\begin{itemize}
\item[(I)] For any recursively enumerable $\mathcal{L}$-theory $T$, there exists $\mathcal{M} \models T$ with $\langle L_\alpha, \in \rangle \prec_{e, n} \mathcal{M}$ and $\mathrm{Ord}^\mathcal{M} \backslash \alpha$ is nonempty and has no least element.
\item[(II)] $\langle L_\alpha, \in \rangle \models \Pi_n\textsf{-Collection}$. 
\end{itemize}
\Square
\end{Coroll1}
   
We suspect that Barwise's admissible cover machinery \cite[Appendix]{bar75} may useful in shedding light on the following question:

\begin{Quest1}
Is there a version of Theorem \ref{Th:MainEndExtensionResult} that holds for all countable models of $\mathsf{KP}+\Pi_n\textsf{-Collection}+\Pi_{n+1}\textsf{-Foundation}$?
\end{Quest1}          

\bibliographystyle{alpha}
\bibliography{.}        

\begin{thebibliography}{9}

\bibitem[Bar75]{bar75} Barwise, J. {\it Admissible Sets and Structures}. Perspectives in Mathematical Logic. Springer-Verlag, Berlin-Heidelberg-New York. 1975.

\bibitem[Fri]{fri73} Friedman, H. M. ``Countable models of set theories". \emph{Cambridge Summer School in Mathematical Logic, August 1--21, 1971}. Edited by A. R. D. Mathias and H. Rogers Jr. Springer Lecture Notes in Mathematics. Vol. 337. Springer, Berlin. 1973. pp 539--573.

\bibitem[FLW]{flw16} Friedman, S.-D., Li, W.; and Wong, T. L. ``Fragments of Kripke-Platek Set Theory and the Metamathematics of $\alpha$-Recursion Theory". \emph{Archive for Mathematical Logic}. Vol. 55. No. 7. 2016. pp 899--924.

%\bibitem[Hut]{hut76} Hutchinson, John E. ``Elementary Extensions of Countable Models of Set Theory". \emph{The Journal of Symbolic Logic}. Vol. 41. No. 1. 1976. pp 139--145.

\bibitem[Gos]{gos80} Gostanian, R. ``Constructible models of subsystems of $\mathrm{ZF}$". {\it The Journal of Symbolic Logic}. Vol. 45. No. 2. 1980. pp 237--250.

\bibitem[JS]{js18} J\"{a}ger, G.; and Steila, S. ``About some fixed point axioms and related principles in Kripke-Platek environments". {\it The Journal of Symbolic Logic}. Vol. 83. No. 2. 2018. pp 642--668.

\bibitem[Kau]{kau81} Kaufmann, M. ``On Existence of $\Sigma_n$ End Extensions". \emph{Logic Year 1979-80, The University of Connecticut}. Lecture Notes in Mathmeatics. No. 859. Springer-Verlag. 1981. pp 92--103.

\bibitem[KM]{km68} Keisler, H. J.; and Morley, M. ``Elementary extensions of models of set theory". {\it Israel Journal of Mathematics}. Vol. 5. 1968. pp 49--65

\bibitem[MS]{ms61} MacDowell, R.; and Specker, E. ``Modelle der Arithmetik" in {\it Infinitistic Methods}. Proceedings of the Symposium on the Foundations of Mathematics, Warsaw, September 1959. Pergamon Press, Oxford and Panstowa Wydanictwo Naukowe, Warsaw. 1961. pp 257--263.

\bibitem[M]{mck19} McKenzie, Z. ``On the relative strengths of fragments of collection". {\it Mathematical Logic Quarterly}. Vol. 65. No. 1. 2019. pp 80--94.

%\bibitem[Mat69]{mat69} Mathias, Adrian R. D. ``Notes on set theory". Available online: {\tt http://www.dpmms.cam.ac.uk/\textasciitilde ardm/} (last accessed on 29/iv/2018)

\bibitem[Mat01]{mat01} Mathias, A. R. D. ``The strength of Mac Lane set theory". \emph{Annals of Pure and Applied Logic}. Vol. 110. 2001. pp 107-234.

\bibitem[PK]{pk78} Paris, J. B.; and Kirby, L. A. S. ``$\Sigma_n$-collection schemas in arithmetic". In {\it Logic Colloquium '77 (Proceedings of the colloquium held in Wroc\l aw, August 1977)}. {\it Studies in Logic and the Foundations of Mathematics}. Vol. 96. North-Holland, Amsterdam-New York, 1978. pp 199--209.

\bibitem[Res]{res87} Ressayre, J.-P. ``Mod\`{e}les non standard et sous-syst\`{e}mes remarquables de ZF". In {\it Mod\`{e}les non standard en arithm\'{e}tique et th\'{e}orie des ensembles}. Volume 22 of {\it Publications Math\'{e}matiques de l'Universit\'{e} Paris VII}. Universit\'{e} de Paris VII, U.E.R. de Math\'{e}matiques, Paris, 1987. pp 47--147. 

\bibitem[Tak]{tak72} Takahashi, M. ``$\tilde{\Delta}_1$-definability in set theory". \emph{Conference in mathematical logic --- London '70}. Edited by W. Hodges. Springer Lecture Notes in Mathematics. Vol. 255. Springer. 1972. pp 281-304.

\end{thebibliography}

\end{document}